\documentclass{article}
\pdfoutput=1

\usepackage{graphicx}
\usepackage{amsthm}
\usepackage{amsmath, amssymb}
\usepackage{hyperref}
\usepackage[title]{appendix}

\newtheorem{thm}{Theorem}
\newtheorem{lem}{Lemma}
\newtheorem{dfn}{Definition}
\newtheorem{ex}{Example}
\newtheorem{rem}{Remark}

\newcommand{\RR}{{\mathcal R}}
\newcommand{\R}{{\mathbb R}}
\newcommand{\ad}{{\mathrm{ad}}}

\begin{document}
\title{Invariants of linear control systems with analytic matrices and the linearizability problem
\thanks{This work was financially supported by Polish National Science Centre grant no. 2017/25/B/ST1/01892.}}

\author{Katerina~V.~Sklyar\thanks{K.V. Sklyar is with Institute of Mathematics, University of Szczecin, Wielkopolska str. 15, Szczecin 70-451, Poland (e-mail: jekatierina.sklyar@usz.edu.pl).}  \ 
		and
        Svetlana~Yu.~Ignatovich\thanks{S.Yu. Ignatovich is with V.N. Karazin Kharkiv National University, Svobody sqr. 4, Kharkiv 61022, Ukraine (e-mail: ignatovich@ukr.net).} 
}

\date{}

\maketitle

\begin{abstract} The paper continues the authors' study of the linearizability problem for nonlinear control systems. In the recent work [K. Sklyar, Systems Control Lett. 134 (2019), 104572], conditions on mappability of a nonlinear control system to a preassigned linear system with analytic matrices were obtained. In the present paper we solve more general problem on linearizability conditions without indicating a target linear system. To this end, we give a description of invariants for linear non-autonomous single-input controllable systems with analytic matrices, which allow classifying such systems up to transformations of coordinates. This study leads to one problem from the theory of linear ordinary differential equations with meromorphic coefficients. As a result, we obtain a criterion for mappability of nonlinear control systems to linear control systems with analytic matrices.
\end{abstract}


\section{Introduction}

The problem of linearizability attracts a great attention of experts in the control theory during almost fifty years. Starting from A.~Krener
\cite{Krener},
in the study of linearizability the Lie bracket technique is applied. We mention several important first publications
\cite{Brockett,Jakubczyk_Respondek,Su_1982,Resp,Respondek_1985};
this list undoubtedly is far from being complete. An important analytic tool here is the Frobenius Theorem on solvability of a system of partial differential equations of the first order. Within this approach, traditionally, control systems of class $C^\infty$ were considered. An alternative idea was introduced by V.~I.~Korobov
\cite{Korobov},
who proposed to consider a special class of triangular systems. This approach admitted $C^1$-smooth systems. Later, these ideas were extended to general nonlinear control systems of class $C^1$
\cite{SCL_2005,CMA_2014,SCL_2016,MCCA_2017,OT_2018}.

However, most researchers focused on autonomous control systems when studying the linearizability. In the present paper we concentrate on a non-autonomous linearizability problem. More specifically, let us consider a nonlinear control system 
\begin{displaymath} 
\dot x=f(t,x,u), \ x\in Q\subset \R^n, \ u\in \R, \ t\in[\alpha,\beta],
\end{displaymath}
where $Q$ is a given domain and $[\alpha,\beta]$ is a given time interval. We are interested in mapping such a system to a linear non-autonomous system
\begin{equation}\label{eq_intr}
\dot z=A(t)z+b(t)u
\end{equation}
by use of a change of variables assuming that it is also time-dependent, i.e., $z=F(t,x)$. Following 
\cite{SCL_2005},
we assume that the vector function $f(t,x,u)$ is of class $C^1$ and a change of variables $F(t,x)$ is of class $C^2$. We also require that a target linear system has analytic matrices, i.e., $A(t)$ and $b(t)$ are real analytic on the interval $[\alpha,\beta]$. Such systems were considered, e.g., in
\cite{KS1,JMAA}.
It was shown that the Markov moment problem 
\cite{M,KN} 
can be efficiently applied for solving the time-optimal control problem for such systems. Moreover, conditions were proposed under which the optimal control can be found by the successive approximation method, on each step of which a power Markov moment problem is solved. In turn, a method of explicit solving the power Markov moment problem was proposed in
\cite{KS};
see also 
\cite{KSI}.

The non-autonomous statement of the linearizability problem mentioned above was proposed in 
\cite{K},
where conditions on mappability of a nonlinear control system to a \emph{preassigned} linear system \eqref{eq_intr} with analytic matrices were obtained. In the present paper we solve more general problem: to find conditions under which a nonlinear system is mapped to \emph{some} linear system with analytic matrices, which is not given in advance.

The precise formulation of the problem is given in Section~\ref{back}. As a preparatory step, in Section~\ref{main} we study some properties of the class of linear non-autonomous systems with analytic matrices, which are closely connected with the classical study of homogeneous linear differential equations of order $n$ with meromorphic coefficients
\cite{Forsyth}. 
In Section~\ref{next} we apply the obtained results and obtain the theorem on non-autonomous linearizability of nonlinear control systems.

\section{Background and statement of the problem}\label{back}

Let us consider a linear control system of the form 
\begin{equation}\label{lin1}
\dot x = A(t)x+b(t)u
\end{equation}
where the matrix $A(t)$ and the vector $b(t)$ are real analytic on the interval $[\alpha,\beta]$. For this system, we introduce the following matrix
\begin{displaymath}
K(t)=(\Delta^0(t),\Delta^1(t),\ldots,\Delta^{n-1}(t)), \  \mbox{where} \  
\Delta^k(t)=(-A(t)+{\textstyle \frac{d}{dt}})^kb(t), \ k\ge 0.
\end{displaymath}
Below we assume that the system \eqref{lin1} is controllable on the time interval $[\alpha,\beta]$. Hence, the matrix $K(t)$ is invertible on $[\alpha,\beta]$ except maybe a finite number of points, say, $\{t_k\}_{k=1}^N$. Let us introduce the following vector function 
\begin{equation}\label{gamma} 
\gamma(t)=K^{-1}(t)\Delta^{n}(t). 
\end{equation}
Obviously, its components are meromorphic functions on $[\alpha,\beta]$.

Suppose we apply a linear change of variables $z=F(t)x$ in the system \eqref{lin1}, where $F(t)$ is a nonsingular analytic matrix, and the system in the new variables takes the form
\begin{displaymath}
\dot z = \widetilde A(t)z+\widetilde b(t)u.
\end{displaymath}
Then $\widetilde A(t)=(\dot F(t)+F(t)A(t))F^{-1}(t)$, $\widetilde b(t)=F(t)b(t)$. Denote
\begin{displaymath}
\widetilde K(t)=(\widetilde\Delta^0(t),\widetilde\Delta^1(t),\ldots,\widetilde\Delta^{n-1}(t)),  \  \mbox{where} \   \widetilde\Delta^k(t)=(-\widetilde A(t)+{\textstyle \frac{d}{dt}})^k\widetilde b(t), \ k\ge 0,
\end{displaymath}
\begin{displaymath}
\widetilde \gamma(t)=\widetilde K^{-1}(t)\widetilde \Delta^{n}(t). 
\end{displaymath}
Then $\widetilde \Delta^k(t)=F(t)\Delta^k(t)$, $k\ge0$. Hence, $\widetilde K(t) =F(t) K(t)$ is invertible on $[\alpha,\beta]\backslash\{t_i\}_{i=1}^N$ and
\begin{displaymath}
\widetilde \gamma(t)=\gamma(t), \ \ t\in[\alpha,\beta]\backslash\{t_i\}_{i=1}^N.
\end{displaymath} 
Therefore, the components of $\gamma(t)$ are invariant w.r.t. (linear analytic) changes of  variables. Following
\cite{K},
we call them \emph{invariants} of the system \eqref{lin1}.

If the initial system is autonomous, $\dot x=Ax+bu$, then $K(t)$ turns into the Kalman matrix $K=(b,-Ab,\ldots,(-A)^{n-1}b)$. Then $\gamma(t)$ is a constant vector,  $\gamma=K^{-1}(-A)^nb$. Therefore, in this case the invariants are constant and, moreover, they \emph{can be arbitrary numbers}. In other words, any collection of numbers $\gamma_1,\ldots,\gamma_n$ uniquely  defines the set of linear systems that are mapped to each other by linear analytic changes of variables. 

Unlike the autonomous systems, in the general case invariants are meromorphic functions, however, \emph{not every meromorphic function is an invariant} for some linear system. Therefore, the following problem arises.

\textbf{Realizability Problem}.  \emph{For a given set of functions $\gamma_1(t),\ldots,\gamma_n(t)$, which are meromorphic on the interval $[\alpha,\beta]$, to determine if they are invariants for some linear control system of the form \eqref{lin1} with analytic matrices.} 

It turns out that this realizability problem plays an important role in linearizability conditions for nonautonomous nonlinear control systems. More specifically, let us consider an affine control system of the form 
\begin{equation}\label{af}
\dot x=a(t,x)+b(t,x)u,
\end{equation}
where $a(t,x), b(t,x) \in C^1([\alpha,\beta]\times Q)$, $Q\subset\R^n$. In order to introduce an analog of invariants for the system \eqref{af},  let us consider the operator  $\RR$  acting as 
\begin{equation}\label{r_fi}
\RR\varphi(t,x)=\varphi_t(t,x)+\varphi_x(t,x)a(t,x)-a_x(t,x)\varphi(t,x)
\end{equation}
for any vector function $\varphi(t,x)\in C^1([\alpha,\beta]\times Q)$. Assuming that all vector functions 
\begin{equation}\label{rkb}
\RR b(t,x), \RR^2 b(t,x), \ldots, \RR^n b(t,x)
\end{equation}
exist, we introduce the matrix
\begin{displaymath}
R(t,x)=(b(t,x),\RR b(t,x),\ldots,\RR^{n-1}b(t,x))
\end{displaymath}
and the vector function 
\begin{equation}\label{gam_tx}
\gamma(t,x)=R^{-1}(t,x)\RR^{n}b(t,x).
\end{equation}

In the particular case of a \emph{driftless} control system of the form
\begin{equation}\label{d_less}
\dot x=g(t,x)u,
\end{equation}
the operator $\RR(t,x)$ reduces to the derivative on $t$, i.e., $\RR^kg(t,x)=g^{(k)}(t,x)$, $k\ge0$ (here and below $g^{(k)}(t,x)$ means the $k$-th derivative on $t$). In 
\cite{K},
driftless systems \eqref{d_less} were considered as a kind of a canonical form for general affine systems \eqref{af}. Under the assumptions $a(t,x)\in C^2([\alpha,\beta]\times Q)$ and $b(t,x)\in C^1([\alpha,\beta]\times Q)$, the system \eqref{af} can be transformed to a dristless form \eqref{d_less} locally, at a neighborhood of any point from $[\alpha,\beta]\times Q$, by a change of variables $z=F(t,x)$ of class $C^2$. A non-local transformation requires additional conditions as well as the case $a(t,x)\in C^1([\alpha,\beta]\times Q)$. In
\cite{K}
driftless systems were used in order to obtain conditions for linearizability of nonautonomous nonlinear systems. 

\begin{dfn}[\cite{K}] \label{def_lin}
We say that a system of the form \eqref{af} is \emph{locally analytically linearizable in the domain $Q$ on the time interval $[\alpha,\beta]$} if there exists a change of variables
\begin{equation}\label{f3}
z=F(t,x)\in C^2([\alpha,\beta]\times Q)
\end{equation}
satisfying the condition 
\begin{equation}\label{f4}
\det F_x(t,x)\neq 0, \ \ (t,x)\in[\alpha,\beta]\times Q,
\end{equation}
such that the system \eqref{af} in the new variables takes the form \eqref{lin1}, where $A(t)$ and $b(t)$ are analytic on $[\alpha,\beta]$.
\end{dfn}

In this definition, the word ``locally'' recalls that the change of variables is only \emph{locally} invertible in the general case and the word ``analytically'' means that a target nonautonomous linear system has analytic matrices. 

\begin{dfn}[\cite{K,K1}] \label{def_lin_pre}
We say that a system of the form \eqref{af} is \emph{locally analytically mappable in the domain $Q$ on the time interval $[\alpha,\beta]$ to a preassigned linear controllable system \eqref{lin1} where $A(t)$, $b(t)$ are analytic on $[\alpha,\beta]$} if there exists a change of variables satisfying \eqref{f3}, \eqref{f4} such that the system \eqref{af} in the new variables takes the form~\eqref{lin1}.
\end{dfn}

We emphasize that, unlike Definition~\ref{def_lin}, in Definition~\ref{def_lin_pre} the target linear system is \emph{fixed}. This requirement is justified by the theorem on linearizability conditions, originally proposed and proved in 
\cite{K}. 
Here we formulate the modification obtained in 
\cite{K1}.

\begin{thm}[\cite{K,K1}] \label{th3}
The system \eqref{af}, where $a(t,x)\in C^2([\alpha,\beta]\times Q)$ and $b(t,x)\in C^1([\alpha,\beta]\times Q)$, is locally analytically mappable in the domain $Q$ on the time interval $[\alpha,\beta]$ to a preassigned linear controllable system \eqref{lin1} if and only if all vector functions \eqref{rkb} exist, belong to the class $C^1([\alpha,\beta]\times Q)$, and satisfy the conditions
\begin{equation}\label{dd1}
[\RR^jb(t,x),\RR^kb(t,x)]=0, \ x\in Q, \ t\in[\alpha,\beta], \ 0\le j,k\le n-1,
\end{equation}
\begin{equation}\label{ddd1}
{\rm rank}\, R(t,x)=n, \ \ t\in [\alpha,\beta]\backslash\{t_i\}_{i=1}^N, \ x\in Q,
\end{equation}
and
\begin{equation}\label{d1}
R^{-1}(t,x)\RR^{n}b(t,x)=K^{-1}(t)\Delta^{n}(t), \ t\in[\alpha,\beta]\backslash\{t_i\}_{i=1}^N.
\end{equation}
\end{thm}

In other words, invariants \eqref{gamma} are unchangeable also under \emph{nonlinear} changes of variables in  linear systems. 

The conditions \eqref{dd1}, \eqref{ddd1} generalize linearizability conditions for autonomous systems. More specifically, if the system \eqref{af} is autonomous, i.e., $a(t,x)=a(x)$ and $b(t,x)=b(x)$, then $\RR^kb(t,x)=\ad^k_{a}b(x)$, $k\ge0$, where $\ad^0_ab(x)=b(x)$, $\ad_a^{k+1}b(x)=[a(x),\ad_a^kb(x)]$, $k\ge0$, and $[\cdot,\cdot]$ means the Lie bracket. The condition  \eqref{d1} arises since the target linear system  \eqref{lin1} is fixed. One may suggest that, dropping the assumption \eqref{d1}, we get the conditions for mappability to a linear system, which is not specified in advance. However, this is not true. Namely, \eqref{d1} should be substituted by the following assumption: components of the vector function $R^{-1}(t,x)\RR^{n}b(t,x)$ depend only on $t$ and \emph{are realizable as invariants of some linear system}. So, we are led to the Realizability Problem mentioned above. 

In the present paper we solve this problem and obtain the conditions of local analytic linearizability (Theorem~\ref{T3}). 

To start with, let us reformulate the Realizability Problem. For a given linear controllable system \eqref{lin1}  with analytic matrices, denote by $\Phi(t)$ the solution of the matrix Cauchy problem $\dot \Phi(t)=A(t)\Phi(t)$, $F(0)=I$. The matrix $\Phi(t)$ is nonsingular and real analytic on $[\alpha,\beta]$. Then in the variables $z=\Phi^{-1}(t)x$ the system \eqref{lin1} takes a driftless form
\begin{equation}\label{lin}
\dot z = \widehat g(t)u,
\end{equation}
where $\widehat g(t)=\Phi^{-1}(t)b(t)$. In this case the calculations mentioned above are simplified since $\Delta^k(t)$ reduces to the $k$-th derivative on $t$. We get 
\begin{displaymath}
\gamma(t)=\widehat K^{-1}(t)\widehat g^{(n)}(t),  \ \ \widehat K(t)=(\widehat g(t),\widehat g^{(1)}(t),\ldots,\widehat g^{(n-1)}(t)),
\end{displaymath}
which can be rewritten as
\begin{equation}\label{syst}
\left(\widehat g(t),\widehat g^{(1)}(t),\ldots,\widehat g^{(n-1)}(t)\right)\gamma(t)=\widehat g^{(n)}(t).  
\end{equation}
This equality means that the column vector $\gamma(t)=(\gamma_1(t),\ldots,\gamma_n(t))^\top$ is a solution of the system of linear \emph{algebraic} equations \eqref{syst}, where the matrix and the right hand side are defined by the known vector function $\widehat g(t)$. 

Now, suppose the contrary: let $\gamma_1(t),\ldots,\gamma_n(t)$ be known and let $\widehat g(t)=(\widehat g_1(t),\ldots,\widehat g_n(t))^\top$ be unknown. Then equalities \eqref{syst} become linear \emph{differential} equations for the components $\widehat g_1(t),\ldots,\widehat g_n(t)$ of the vector function $\widehat g(t)$. More specifically, let us consider the following differential equation
\begin{equation}\label{dif_eq}
\gamma_1(t)y+\gamma_2(t)y^{(1)}+\cdots+\gamma_n(t)y^{(n-1)}=y^{(n)}.  
\end{equation}
Then \eqref{syst} means that this differential equation has $n$ linearly independent real analytic solutions $\widehat g_1(t),\ldots,\widehat g_n(t)$. Therefore, our Realizability Problem can be reformulated in purely classical terms.

\textbf{Realizability Problem (reformulated).} \emph{For a given set of functions $\gamma_1(t),\ldots,\gamma_n(t)$, which are meromorphic on the interval $[\alpha,\beta]$, to determine if the differential equation \eqref{dif_eq} has $n$ linearly independent real analytic solutions on $[\alpha,\beta]$.}

Linear differential equations with meromorphic coefficients were studied in detail for the case $n=2$ due to their significance for the mathematical physics
\cite{Whit_Wat,Codd_Lev,Teschl}.
Here the method of F.G. Frobenius 
\cite{Frobenius}
is applicable: try to find $n$ independent solutions in the form $y(t)=t^{a_i}\varphi_i(t)$, where $a_i$ are  constants and $\varphi_i(t)$ are analytic functions. If this is impossible, the function $\log t$ should be involved. For the case $n>2$, the solution becomes complicated; it is discussed in 
\cite{Forsyth}. 

The Realizability Problem formulated above is different: for the equation \eqref{dif_eq}, in the general $n$-dimensional case, to find conditions under which there exist $n$ linearly independent analytic solutions. We give a solution of this problem in Section~\ref{main} and then apply it to  study linearizability conditions for nonlinear control systems in Section~\ref{next}.  

\section{Analytic solvability of linear differential equations}\label{main}

To simplify the notation, in this section we denote $p_s(t)=-\gamma_{n-s+1}(t)$, $s=1,\ldots,n$, i.e., we write the equation \eqref{dif_eq} as
\begin{equation}\label{eq_n}
y^{(n)}+\sum_{s=1}^{n}p_{s}(t)y^{(n-s)}=0.
\end{equation}
Our nearest goal is to formulate conditions under which the equation \eqref{eq_n} has $n$ linearly independent analytic solutions in a neighborhood of a given point $t=t_0$; without loss of generality we assume $t_0=0$. The following lemma describes necessary conditions for the coefficients $p_s(t)$. 

\begin{lem}\label{lem2}
Suppose the equation \eqref{eq_n} has $n$ linearly independent analytic solutions in a neighborhood of the point $t=0$. Then each $p_s(t)$ is analytic or meromorphic with a pole of order no greater than $s$ at the point $t=0$. 
\end{lem}

As is well known, the necessary condition mentrioned in Lemma~\ref{lem2} follows from much more general requirements
\cite[Ch. IV, Theorem 5.2]{Codd_Lev}. However, in our case the proof is easy; we give it for the sake of completeness.

\emph{Proof.} First, we note that if $y(t)$ is analytic at the point $t=0$, then $\frac{y^{(k)}(t)}{y(t)}$ is analytic or meromorphic with a pole of order at most $k$ at $t=0$.

We argue by induction on $n$. For $n=1$, the lemma is obvious. Now, for $n\ge2$, let us suppose that $y_1(t),\ldots,y_{n}(t)$ are linearly independent analytic solutions of the equation \eqref{eq_n}. Without loss of generality we assume that $y_k(t)=z_k(t)y_{n}(t)$, where $z_k(t)$ are analytic functions, $k=1,\ldots,n-1$. Substituting $y_1(t),\ldots,y_{n-1}(t)$ to the equation \eqref{eq_n},  after obvious simplification we obtain that $\dot z_1(t),\ldots,\dot z_{n-1}(t)$ are $n-1$ linearly independent analytic solutions of the equation 
\begin{displaymath}
y^{(n-1)}+\sum_{i=1}^{n-1}\widetilde p_{i}(t)y^{(n-1-i)}=0,
\end{displaymath}
with
\begin{equation}\label{eq_n+1_00}
\widetilde p_{i}(t)=p_{i}(t)+\sum_{s=1}^{i-1}C_{n-s}^{i-s}\,p_s(t)\frac{y_{n}^{(i-s)}(t)}{y_{n}(t)}+C_n^i\frac{y_{n}^{(i)}(t)}{y_{n}(t)}, \ \ i=1,\ldots,n-1,
\end{equation}
where $C_j^i=\frac{j!}{i!(j-i)!}$. By the induction supposition, every $\widetilde p_{i}(t)$ has a pole up to order $i$ at $t=0$. Then, using the induction on $i$, we obtain from \eqref{eq_n+1_00} that every $p_i(t)$ also has a pole of order at most $i$ at $t=0$ for all $i=1,\ldots,n-1$. Finally, substituting $y_n(t)$ to \eqref{eq_n} and expressing $p_{n}(t)$ as
\begin{displaymath}
p_{n}(t)=-\frac{y_{n}^{(n)}(t)}{y_{n}(t)}-\sum_{s=1}^{n-1}p_s(t)\frac{y_{n}^{(n-s)}(t)}{y_{n}(t)},
\end{displaymath}
we get that $p_{n}(t)$ has a pole of order at most $n$ at $t=0$. 
\qed

Thus, we restrict ourselves by those $p_s(t)$ that have poles up to order $s$. Suppose $y(t)$ is an analytic solution of the equation \eqref{eq_n}. Let us expand $p_s(t)$ and $y(t)$ into series on $t$, 
\begin{displaymath}
p_s(t)=\sum_{i=-s}^\infty p_{s,i}t^i, \ \ s=1,\ldots,n, \ \ 
y(t)=\sum_{k=0}^\infty y_kt^k.
\end{displaymath}

Below we use the notation for \emph{falling factorial} for nonnegative integers $k,q$
\cite[Subsection 2.6]{Knuth},
\begin{displaymath}
k^{\underline{q}}=\left\{\begin{array}{cl}
\displaystyle \frac{k!}{(k-q)!} &\mbox{if} \ q\le k,\\
\displaystyle 0 &\mbox{if} \  q>k.
\end{array} \right.
\end{displaymath}
In other words, for any integer $k\ge0$, 
\begin{displaymath}
k^{\underline{0}}=1 \ \mbox{ and } \ k^{\underline{q}}=\underbrace{k(k-1)\cdots(k-q+1)}_{q \ \mbox{\footnotesize multipliers}} \ \mbox{for} \ q\ge1.
\end{displaymath}
Then 
\begin{displaymath}
y^{(q)}(t)=\sum_{k=0}^\infty k^{\underline{q}} y_kt^{k-q}.
\end{displaymath}

Substituting the series for $p_s(t)$ and $y^{(n-s)}(t)$ to \eqref{eq_n}, we get
\begin{displaymath}
\sum_{k=0}^\infty k^{\underline{n}}y_kt^{k-n}+\sum_{s=1}^{n}\sum_{i=-s}^\infty p_{s,i}t^i\sum_{j=0}^\infty j^{\underline{n-s}}y_jt^{j-n+s}=0.
\end{displaymath}
Finding the coefficients of powers of $t$ in the left hand side and equating them to zero, we get the system of equations
\begin{equation}\label{eq_s1}
k^{\underline{n}}y_k+\sum_{j=0}^{k}\sum_{s=1}^{n}j^{\underline{n-s}}p_{s,k-j-s}y_j=0, \ \ k\ge 0.
\end{equation}
Let us introduce  the notation 
\begin{equation}\label{rk_gen}
\begin{array}{l}
\displaystyle
V_{k,k}=k^{\underline{n}}+\sum_{s=1}^{n}k^{\underline{n-s}}p_{s,-s} \ \ \mbox{ for} \ \ k\ge0, \\ 
\displaystyle
V_{k,j}=\sum_{s=1}^{n}j^{\underline{n-s}}p_{s,k-j-s} \ \ \mbox{ for} \ \ k\ge0, \ \ 0\le j\le k-1,\\ 
\end{array}
\end{equation}
then the system  \eqref{eq_s1} takes the form 
\begin{equation}\label{eq_s1_new}
V_{k,k}y_k+\sum_{j=0}^{k-1}V_{k,j}y_j=0, \ \ k\ge 0.
\end{equation}
If $V_{k,k}\not=0$, then $y_k$ is uniquely defined by $y_0,\ldots,y_{k-1}$. The equality $V_{k,k}=0$, called the \emph{indicial equation}, has the form 
\begin{equation}\label{eq_k}
k^{\underline{n}}+\sum_{s=1}^{n}k^{\underline{n-s}}p_{s,-s}=0;
\end{equation}
it is a polynomial equation for $k$ of degree $n$. Therefore, if the system \eqref{eq_s1_new} has $n$ linearly independent solutions, then the equation \eqref{eq_k} has $n$ different nonnegative integer roots. 

Suppose these roots are $0\le k_1<\cdots<k_n$. Then equations \eqref{eq_s1_new} for $k\le k_1-1$ give $y_0=\cdots=y_{k_1-1}=0$ while equations \eqref{eq_s1_new} for $k\ge k_n+1$ uniquely define all $y_k$ for $k\ge k_n+1$. Since $y_0=\cdots=y_{k_1-1}=0$ and $V_{k_1,k_1}=0$, the equation  \eqref{eq_s1_new} for $k=k_1$ is trivial. Let us write the rest equations \eqref{eq_s1_new} for $k=k_1+1,\ldots,k_n$. Taking into account that $V_{k_n,k_n}=0$ and using the matrix notation, we have
\begin{equation}\label{det=0_gen_syst}
\begin{pmatrix}
V_{k_1+1,k_1}&V_{k_1+1,k_1+1}&0&\cdots&0\\
\cdots&\cdots&\cdots&\cdots&\cdots\\
V_{k_{n}-1,k_1\!\!}&V_{k_{n}-1,k_1+1\!\!}&\cdots&\cdots&V_{k_{n}-1,k_{n}-1\!\!}\\
V_{k_{n},k_1}&V_{k_{n},k_1+1}&\cdots&\cdots&V_{k_{n},k_{n}-1}\\
\end{pmatrix}
\begin{pmatrix}
y_{k_1}\\\cdots\\y_{k_{n}-2}\\y_{k_{n}-1}
\end{pmatrix}=0.
\end{equation}
This is the system of $k_n-k_1$ linear equations in $k_n-k_1$ unknowns. We note that this system does not include $y_{k_n}$, which can be arbitrary. Therefore, we are interested in condictions under which the system \eqref{det=0_gen_syst} has $n-1$ linearly independent solutions. Obviously, this is the case  if and only if 
\begin{equation}\label{det=0_gen}
{\rm rank}\begin{pmatrix}
V_{k_1+1,k_1}&V_{k_1+1,k_1+1}&0&\cdots&0\\
\cdots&\cdots&\cdots&\cdots&\cdots\\
V_{k_{n}-1,k_1\!\!}&V_{k_{n}-1,k_1+1\!\!}&\cdots&\cdots&V_{k_{n}-1,k_{n}-1\!\!}\\
V_{k_{n},k_1}&V_{k_{n},k_1+1}&\cdots&\cdots&V_{k_{n},k_{n}-1}\\
\end{pmatrix}=k_n-k_1-n+1.
\end{equation}
For $n=2$, this rank is $k_2-k_1-1$, i.e., equals the dimension of the matrix minus~1. Hence, the condition \eqref{det=0_gen} for $n=2$ holds if and only if the matrix in \eqref{det=0_gen} is singular. In the general case the condition \eqref{det=0_gen} reduces to $\frac{n(n-1)}{2}$ equalities, see Remark~\ref{rem1}. For $n=1$, such a condition is not required.

We are led to the following result.

\begin{thm}\label{T2} The equation \eqref{eq_n} has $n$ linearly independent analytic solutions in a neighborhood of the point $t=0$ if and only if the following conditions are satisfied: the polynomial equation \eqref{eq_k} has $n$ different nonnegative integer roots and the condition \eqref{det=0_gen} holds, where $V_{k,j}$ are defined by \eqref{rk_gen}. If this is the case, the components $y_{k_1},\ldots,y_{k_n}$ can be chosen arbitrary and then all other components are defined uniquely by \eqref{eq_s1_new}.
\end{thm}

\emph{Proof.}
\emph{Necessity} is shown above. \emph{Sufficiency}. Suppose that $k_1<\cdots<k_n$ are different nonnegative integer roots of the equation \eqref{eq_k} and the condition \eqref{det=0_gen} holds. Since $V_{k,k}\not=0$ for $k\not=k_s$, $s=2,\ldots,n-1$, exactly $k_n-k_1-n+1$ elements in the first super-diagonal of the matrix in \eqref{det=0_gen} are nonzero. Obviously, the rows containing these elements are linearly independent. Therefore, the condition \eqref{det=0_gen} implies that all other rows, namely, $(V_{k_s,k_1},\ldots,V_{k_s,k_s-1},0,\ldots,0)$ for $s=2,\ldots,n$, linearly depend on previous ones. Thus, $y_{k_1},\ldots,y_{k_n}$ can be chosen arbitrarily while all the rest $y_k$ are defined uniquely from \eqref{eq_s1_new} by the equalities 
\begin{equation}\label{yk_2_gen}
y_k=-\frac{1}{V_{k,k}}\sum_{j=0}^{k-1}V_{k,j}y_j, \ \ k\not=k_1,\ldots,k\not=k_n.
\end{equation}

Let us prove that for any choice of $y_{k_1}, \ldots,y_{k_n}$ the series $\sum_{k=0}^\infty y_kt^k$ converges in a neighborhood of the point $t=0$. Due to our assumptions, there exist $C\ge1$ such that 
\begin{displaymath}
|p_{s,i}|\le C^{s+i}, \ s=1,\ldots,n, \ \ i\ge -s+1
\end{displaymath}
($p_{s,-s}$ are not included into formulas for $V_{k,j}$, $0\le j\le k-1$, see \eqref{rk_gen}). Then
\begin{equation}\label{yk_0_gen_00}
|V_{k,j}|\le \sum\limits_{s=1}^{n}j^{\underline{n-s}}|p_{s,k-j-s}|\le \sum\limits_{s=1}^{n}j^{\underline{n-s}}C^{k-j}.
\end{equation}
Below we use the following identity 
\begin{equation}\label{yk_00_gen}
\sum_{j=0}^{k-1}j^{\underline{m}}=\frac{k^{\underline{m+1}}}{m+1}, 
\end{equation}
where $m\ge0$ and $k\ge1$ are integers; it can be proved easily by induction on $k$ 
\cite[Subsection 2.6]{Knuth}. 

Let us introduce the polynomial 
\begin{equation}\label{yk_30_gen}
P(k)=\sum\limits_{s=1}^{n}\frac{k^{\underline{n-s+1}}}{n-s+1}.
\end{equation}
It is of degree $n$ and its leading coefficient equals $\frac{1}{n}$. Hence, 
\begin{displaymath}
\frac{P(k)}{V_{k,k}}=\frac{P(k)}{\prod\limits_{i=1}^n(k-k_i)}\to\frac1n \ \mbox{ as } \ k\to\infty.
\end{displaymath}

Now suppose $n\ge2$ (the case $n=1$ is considered below). Then there exists $k_0\ge k_n$ such that 
\begin{equation}\label{yk_3_gen}
\frac{P(k)}{V_{k,k}}\le1 \ \  \mbox{for} \ k\ge k_0.
\end{equation}
Let us introduce $C_1$ as
\begin{displaymath}
C_1=\max_{0\le k\le k_0}|y_k|,
\end{displaymath}
then, since $C\ge1$, we get 
\begin{equation}\label{k0}
|y_k|\le C_1C^{k} \ \  \mbox{for} \ 0\le k\le k_0.
\end{equation}
We prove that $|y_k|\le C_1C^k$ for all $k\ge k_0$, by induction on $k$. For $k=k_0$, we have \eqref{k0}. Suppose $k\ge k_0+1$ and 
\begin{equation}\label{yk_4}
|y_j|\le C_1C^{j} \ \ \mbox{for all} \ 0\le j\le k-1.
\end{equation}
Then using \eqref{yk_2_gen}--\eqref{yk_3_gen} and \eqref{yk_4}, we get
\begin{displaymath}
|y_k|\le \frac{\sum\limits_{j=0}^{k-1}\sum\limits_{s=1}^{n}j^{\underline{n-s}}C^{k-j}C_1C^{j}}{V_{k,k}}=C_1C^{k}\frac{\sum\limits_{s=1}^{n}\left(\sum\limits_{j=0}^{k-1}j^{\underline{n-s}}\right)}{V_{k,k}}=
\frac{C_1C^{k}P(k)}{V_{k,k}}\le C_1C^{k}.
\end{displaymath}

In the case $n=1$ we have $P(k)=k$, therefore, \eqref{yk_3_gen} does not hold if $k_1>0$. However, let us take into account that $y_0=\cdots=y_{k_1-1}=0$. We choose $C_1=|y_{k_1}|$ and use the induction supposition \eqref{yk_4} with $k_0=k_1$ and the inequality \eqref{yk_0_gen_00}, which gives $|V_{k,j}|\le C^{k-j}$. Then we directly obtain from \eqref{yk_2_gen}
\begin{displaymath}
|y_k|\le \frac{\sum\limits_{j=0}^{k-1}C^{k-j}|y_j|}{V_{k,k}}\le \frac{\sum\limits_{j=k_1}^{k-1}C^{k-j}C_1C^j}{k-k_1}= C_1C^{k}.
\end{displaymath}

Thus, we proved by induction that $|y_k|\le C_1C^{k}$ for all $k\ge0$, therefore, the series $\sum_{k=0}^\infty y_kt^k$ converges if $|t|<C^{-1}$. 

Therefore, choosing $n$ linearly independent tuples $(y_{k_1},\ldots,y_{k_n})$, we get $n$ linearly independent analytic solutions of the equation \eqref{eq_n}. 
\qed

\begin{ex}\label{ex3}\rm
For $n=2$, the indicial equation \eqref{eq_k} has the form 
\begin{equation}\label{eq_k_n=2}
k(k-1)+kp_{1,-1}+p_{2,-2}=0
\end{equation}
Suppose $p_{1,-1}=-1$ and \eqref{eq_k_n=2} has two nonnegative integer roots $k_1<k_2$. Then $k_1+k_2=-(p_{1,-1}-1)=2$, therefore, the unique possible case is $k_1=0$, $k_2=2$. Then $p_{2,-2}=k_1k_2=0$.

Then \eqref{rk_gen} implies $V_{1,0}=p_{2,-1}$, $V_{1,1}=p_{1,-1}+p_{2,-2}=-1$, $V_{2,0}=p_{2,0}$, $V_{2,1}=p_{1,0}+p_{2,-1}$, therefore, the condition  \eqref{det=0_gen} takes the form
\begin{displaymath}
{\rm det}\begin{pmatrix}
V_{1,0}&V_{1,1}\\
V_{2,0}&V_{2,1}\\
\end{pmatrix}={\rm det}\begin{pmatrix}
p_{2,-1}&-1\\
p_{2,0}&p_{1,0}+p_{2,-1}\\
\end{pmatrix}=0,
\end{displaymath}
which gives $p_{2,-1}(p_{1,0}+p_{2,-1})+p_{2,0}=0$. We note that for any $p_{1,0}$ and $p_{2,-1}$ there exists a unique $p_{2,0}$ satisfying this equality. 

For example, if $p_{1,0}=0$ and $p_{2,-1}=1$, then $p_{2,0}=-1$. If $p_{1,k}=p_{2,k}=0$ for all $k\ge1$, i.e., $p_1(t)=-\frac1t$ and $p_2(t)=-1+\frac1t$, the equation \eqref{eq_n} takes the form
\begin{equation}\label{ex2_eq}
\ddot y-\frac{1}{t}\dot y-\left(1-\frac{1}{t}\right)y=0.
\end{equation}
Two linearly independent analytic solutions can be found from the formula \eqref{yk_2_gen}, which implies
\begin{equation}\label{ex3_eq}
y_1=y_0, \ \ \ y_k=\frac{y_{k-2}-y_{k-1}}{k(k-2)}, \ k\ge3,
\end{equation}
and $y_0$, $y_2$ can be chosen arbitrarily. In this case $C=1$, therefore, the obtained analytic solutions exist at least for $|t|<1$. It is easy to check that the sequences
\begin{displaymath}
y_k=\frac{1}{k!}, \ k\ge0,  \ \mbox{and} \ \
y_k=\frac{(-1)^k(1-2k)}{k!}, \ k\ge0,
\end{displaymath}
satisfy the equation \eqref{ex3_eq}; they are coefficients of the series for $y(t)=e^t$ and $y(t)=e^{-t}(2t+1)$, which are analytic linearly independent solutions of \eqref{ex2_eq}.
\end{ex}

\begin{ex}\label{ex6}\rm
For $n=3$, the indicial equation \eqref{eq_k} has the form 
\begin{displaymath}
k(k-1)(k-2)+k(k-1)p_{1,-1}+kp_{2,-2}+p_{3,-3}=0.
\end{displaymath}
Suppose this equation has three roots $k_1=1$, $k_2=2$, $k_3=4$, which is true if $p_{1,-1}=-4$, $p_{2,-2}=8$, $p_{3,-3}=-8$. In this case condition \eqref{det=0_gen} reads 
\begin{displaymath}
{\rm rank}\begin{pmatrix}
V_{2,1}&0&0\\
V_{3,1}&V_{3,2}&-2\\
V_{4,1}&V_{4,2}&V_{4,3}
\end{pmatrix}=1,
\end{displaymath}
which holds if and only if
\begin{equation}\label{det=0_gen_ex}
V_{2,1}=0, \ \ {\rm det}\begin{pmatrix}
V_{3,2}&-2\\
V_{4,2}&V_{4,3}
\end{pmatrix}=0, \ \ {\rm det}\begin{pmatrix}
V_{3,1}&-2\\
V_{4,1}&V_{4,3}
\end{pmatrix}=0.
\end{equation}
\end{ex}

\begin{rem}\label{rem1}\rm
As was mentioned above, the condition \eqref{det=0_gen} reduces to $\frac{n(n-1)}{2}$ equalities. It is useful to express them as conditions on minors of the matrix from  \eqref{det=0_gen} analogously to the conditions \eqref{det=0_gen_ex} in Example~\ref{ex6}. In order to formulate them, let us denote by $D_{i,j}$ the determinant of the matrix formed by deleting the rows and columns containing $V_{k_s,k_s}$ for all $s$ such that $i<s<j$ from the matrix
\begin{displaymath}
\begin{pmatrix}
V_{k_i+1,k_i}&V_{k_i+1,k_i+1}&0&\cdots&0\\
\cdots&\cdots&\cdots&\cdots&\cdots\\
V_{k_{j}-1,k_i\!\!}&V_{k_{j}-1,k_i+1\!\!}&\cdots&\cdots&V_{k_{j}-1,k_{j}-1\!\!}\\
V_{k_{j},k_i}&V_{k_{j},k_i+1}&\cdots&\cdots&V_{k_{j},k_{j}-1}
\end{pmatrix}
\end{displaymath}
(since $j-i-1$ rows and columns should be deleted, such a matrix is of dimension $k_j-k_i-(j-i-1)$). One can show that \emph{the condition \eqref{det=0_gen} holds if and only if }
\begin{displaymath}
D_{i,j}=0 \ \mbox{ for all } \ 1\le i<j\le n.
\end{displaymath}
So, in Example~\ref{ex6}, $D_{1,2}$, $D_{2,3}$, and $D_{1,3}$ should vanish, which coincides with \eqref{det=0_gen_ex}.
\end{rem}

\begin{rem}\label{rem2}\rm
Let us consider the case when all functions $p_1(t),\ldots,p_n(t)$ are analytic, i.e., $p_{s,j}=0$ for all $-s\le j\le-1$, $s=1,\ldots,n$. Then the indicial equation \eqref{eq_k} takes the form $k^{\underline{n}}=0$; its solutions are $k_i=i-1$, $i=1,\ldots,n$. Taking into account that $j^{\underline{n-s}}=0$ if $n-s>j$, we conclude from \eqref{rk_gen} that  $V_{k,j}=0$ for all $0\le j\le k-1 \le n-2$. Hence, the matrix in \eqref{det=0_gen} is zero. Thus, in this case all the conditions of Theorem~\ref{T2} are trivially satisfied.
\end{rem}

\section{Conditions for local analytic linearizability}\label{next}
Now we return to linear and nonlinear control systems. First, as a corollary of Theorem~\ref{T2}, we obtain a solution of the Realizability Problem. 

\begin{thm}[On realizability]\label{th_r}
Let the functions $\gamma_1(t),\ldots,\gamma_n(t)$ be analytic or meromorphic on the interval $[\alpha,\beta]$. Denote by $\{t_i\}_{i=1}^N\subset [\alpha,\beta]$ the set of points where at least one of them has a pole. These functions are invariants for some linear control system of the form \eqref{lin1} with analytic matrices on $[\alpha,\beta]$ if and only if they satisfy the following conditions at any $t_i$, $i=1,\ldots,N$:

(i) each function $\gamma_s(t)$ is analytic in a neighborhood of $t=t_i$ or meromorphic with a pole at $t=t_i$ of order no greater than $n-s+1$, i.e., $\gamma_s(t)$ are expanded into convergent series
\begin{displaymath}
\gamma_s(t)=\sum_{j=-n+s-1}^\infty \gamma_{s,j}(t-t_i)^j, \ \ s=1,\ldots, n;
\end{displaymath}

(ii) the polynomial equation 
\begin{equation}\label{eq_th}
k^{\underline{n}}-\sum_{s=1}^{n}k^{\underline{n-s}}\gamma_{n-s+1,-s}=0
\end{equation}
has $n$  different nonnegative integer solutions $0\le k_1<\cdots<k_n$; 

(iii) the following equality holds
\begin{equation}\label{eq1_th}
{\rm rank}\begin{pmatrix}
V_{k_1+1,k_1}&V_{k_1+1,k_1+1}&0&\cdots&0\\
\cdots&\cdots&\cdots&\cdots&\cdots\\
V_{k_{n}-1,k_1\!\!}&V_{k_{n}-1,k_1+1\!\!}&\cdots&\cdots&V_{k_{n}-1,k_{n}-1\!\!}\\
V_{k_{n},k_1}&V_{k_{n},k_1+1}&\cdots&\cdots&V_{k_{n},k_{n}-1}\\
\end{pmatrix}=k_n-k_1-n+1,
\end{equation}
where
\begin{displaymath}
\begin{array}{l}
\displaystyle
V_{k,k}=k^{\underline{n}}-\sum_{s=1}^{n}k^{\underline{n-s}}\gamma_{n-s+1,-s} \ \ \mbox{ for} \ \  k_1+1\le k\le k_n-1, \\ 
\displaystyle
V_{k,j}=-\sum_{s=1}^{n}j^{\underline{n-s}}\gamma_{n-s+1,k-j-s} \ \ \mbox{ for} \ \ k_1\le j\le k-1\le k_n-1. 
\end{array}
\end{displaymath}
\end{thm}

\emph{Proof.}
\emph{Necessity} follows from Theorem~\ref{T2}. To prove \emph{sufficiency}, let us consider the equation \eqref{eq_n} with $p_s(t)=-\gamma_{n-s+1}(t)$, $t\in [\alpha,\beta]$, $s=1,\ldots,n$. By our supposition, $p_1(t),\ldots,p_n(t)$ satisfy the conditions of Theorem~\ref{T2} at any point $t\in[\alpha,\beta]$. Hence, any point $t\in[\alpha,\beta]$ has a neighborhood $U(t)$ such that the equation \eqref{eq_n} has $n$ linearly independent analytic solutions in $U(t)$. Let us choose a finite number of points $t_1<\cdots<t_m$ such that $[\alpha,\beta]\subset\bigcup_{r=1}^m U(t_r)$. Without loss of generality we assume that $U(t_r)\cap U(t_{r+1})\not=\varnothing$, $r=1,\ldots,m-1$.

Suppose that $y_1(t),\ldots,y_n(t)$ are $n$ linearly independent analytic solutions of \eqref{eq_n} in $U(t_1)$ and $\widetilde y_1(t),\ldots,\widetilde y_n(t)$ are $n$ linearly independent analytic solutions of \eqref{eq_n} in $U(t_2)$. Let us consider the interval $U(t_1)\cap U(t_2)\not=\varnothing$. Then $y_1(t),\ldots,y_n(t)$ and $\widetilde y_1(t),\ldots,\widetilde y_n(t)$ are two sets of linearly independent solutions of the (linear) differential equation \eqref{eq_n} in $U(t_1)\cap U(t_2)$. Therefore, $y_i(t)=\sum_{j=1}^n a_{ij}\widetilde y_j(t)$, $i=1,\ldots,n$, for  $t\in U(t_1)\cap U(t_2)$, where $a_{ij}$ are some constants and the matrix $\{a_{ij}\}_{i,j=1}^n$ is nonsingular. 

Let us extend the functions $y_i(t)$ to the interval $U(t_2)\backslash U(t_1)$ defining $y_i(t)=\sum_{j=1}^n a_{ij}\widetilde y_j(t)$ for $t\in U(t_2)\backslash U(t_1)$ for any $i=1,\ldots,n$. Then $y_1(t),\ldots,y_n(t)$  become $n$ linearly independent analytic solutions of the equation \eqref{eq_n} in the interval $U(t_1)\cup U(t_2)$. Continuing this process, after a finite  number of steps we obtain $n$ linearly independent analytic solutions of \eqref{eq_n} in  $\bigcup_{r=1}^m U(t_r)$. These solutions, considered as components of the vector function $\widehat g(t)$, generate a linear control system \eqref{lin} with invariants $\gamma_1(t),\ldots,\gamma_n(t)$. 
\qed

\begin{ex}\label{ex_hyp}\rm Consider the functions
\begin{equation}\label{gam_ex_2}
\gamma_1(t)=-\frac{c}{t(1-t)}, \ \ \gamma_2(t)=-\frac{a+bt}{t(1-t)}, \ \ t\in[0,1],
\end{equation}
having poles at $t=0$ and $t=1$. For the point $t=0$,
\begin{displaymath}
\gamma_1(t)=-\frac{c}{t(1-t)}=-\sum_{k=-1}^\infty ct^k, \ \ \gamma_2(t)=-\frac{a+bt}{t(1-t)}=-\frac{a}{t}-\sum_{k=0}^\infty (a+b)t^k.
\end{displaymath}
Therefore, the equation \eqref{eq_th} takes the form $k(k-1)+ka=0$, hence, $k_1=0$, $k_2=1-a$. Thus, $a$ should be a nonpositive integer. One easily find $V_{k,j}=j(a+b)+c$, $0\le j\le k-1$. The condition \eqref{eq1_th} reduces to
\begin{displaymath}
\prod_{j=0}^{-a}(jb+c-j(j-1))=0.
\end{displaymath}
Therefore, an integer $0\le j \le -a$ should exist such that $jb+c=j(j-1)$.

For the point $t=1$, arguing analogously, we get that $a+b$ should be a nonnegative integer and an integer $0\le j \le a+b$ should exist such that $jb+c=j(j-1)$.

As a result, the functions \eqref{gam_ex_2} satisfy the conditions of Theorem~\ref{th_r} if and only if $a\le0$ and $a+b\ge0$ are integers and there exists an integer  $0\le j \le \min\{-a,a+b\}$ such that $jb+c=j(j-1)$.   For example, let $a=-1$, $b=2$, $c=-2$, then $j=1=-a=a+b$ satisfies the condition mentioned above. Therefore, the functions $\gamma_1(t)=\frac{2}{t(1-t)}$ and $\gamma_2(t)=\frac{1-2t}{t(1-t)}$ are invariants for some linear control system of the form \eqref{lin1} with analytic matrices on $t\in [0,1]$. Namely, in this case the equation \eqref{dif_eq}, which is hypergeometric, takes the form
\begin{displaymath}
\frac{2}{t(1-t)}y+\frac{1-2t}{t(1-t)}\dot y=\ddot y.
\noindent
\end{displaymath}
It has two linearly independent solutions, which can be chosen as $y_1(t)=2t-1$ and $y_2(t)=t^2$. Therefore, the functions $\gamma_1(t)$ and $\gamma_2(t)$ are invariants of the system \eqref{lin} with $\widehat g(t)=(2t-1,t^2)^\top$.
\end{ex}

\begin{rem}\rm We emphasize that Theorem~\ref{th_r} allows us to answer the realizability question \emph{without solving the equation \eqref{dif_eq}}. However, suppose that we find $n$ linearly independent solutions of the equation \eqref{dif_eq} that turn out to be analytic in an interval including one or several points from the set $\{t_i\}_{i=1}^N$. Then the conditions of Theorem~\ref{th_r} for all these points are satisfied automatically, hence, we do not need to check them. 
\end{rem}

\begin{ex}\rm
Consider the functions
\begin{equation}\label{eq_no}
\gamma_1(t)=-\frac{2}{t^2(1-t)}, \ \ \gamma_2(t)=\frac{2}{t(1-t)}, \ \  t\in[0,1],
\end{equation}
having poles at $t=0$ and $t=1$. At the point $t=0$ we get $\gamma_{1,-2}=-2$, $\gamma_{2,-1}=2$, therefore, $k_1=1$, $k_2=2$. The condition \eqref{eq1_th} takes the form $V_{2,1}=0$; obviously it holds since $V_{2,1}=-\gamma_{2,0}-\gamma_{1,-1}=-2+2=0$. However, at the point $t=1$ we have $\gamma_{1,-2}=0$, $\gamma_{2,-1}=-2$. Hence, the indicial equation does not have two nonnegative roots. Therefore, the functions \eqref{eq_no} are invariants for a linear control system of the form \eqref{lin1} with analytic matrices on any interval $[0,\beta]$ such that $\beta<1$ but not on the interval $[0,1]$. In order to find such a system, let us 
consider the equation \eqref{dif_eq}, which takes the form
\begin{displaymath}
-\frac{2}{t^2(1-t)}y+\frac{2}{t(1-t)}\dot y=\ddot y.
\end{displaymath}
One can verify that $y_1(t)=t$ and $y_2(t)=\frac{t}{1-t}$ are two linearly independent solutions. Hence, the functions \eqref{eq_no} are invariants of the system \eqref{lin} with $\widehat g(t)=(t,\frac{t}{1-t})^\top$ defined on $[0,\beta]$ with $\beta<1$. 
\end{ex}

As a consequence of Theorems~\ref{th3} and~\ref{th_r}, we obtain conditions for local analytic linearizability.

\begin{thm}[On local analytic linearizability]\label{T3} 
Consider a nonlinear control system of the form \eqref{af}, where  $a(t,x)\in C^2([\alpha,\beta]\times Q)$, $b(t,x)\in C^1([\alpha,\beta]\times Q)$. This system is locally analytically linearizable in the domain $Q$ on the time interval $[\alpha,\beta]$ if and only if all vector functions \eqref{rkb} exist, belong to the class $C^1([\alpha,\beta]\times Q)$, satisfy the conditions \eqref{dd1} and  \eqref{ddd1}, and components of the vector function \eqref{gam_tx} depend only on $t$, i.e., $\gamma(t,x)=\gamma(t)$, and are invariants for some linear control system of the form \eqref{lin1} with analytic matrices on $[\alpha,\beta]$, i.e., satisfy the conditions of Theorem~\ref{th_r}. 
\end{thm}

\begin{ex}\rm Let us consider the following nonlinear control system
\begin{equation}\label{sys_ex_last}
\dot x_1=\frac{4t^2|t|}{2+t^3|t|}x_1+(2+t^3|t|)(1+3t^3x_2^2)u, \ \ \ 
\dot x_2=t^3u
\end{equation}
of the class $C^2([-1,1]\times \R^2)$. We have
\begin{displaymath}
b(t,x)=\begin{pmatrix}
\rho(t)(1+3t^3x_2^2)\\
t^3
\end{pmatrix},   
\RR b(t,x)=\begin{pmatrix}
9t^2\rho(t)x_2^2\\
3t^2
\end{pmatrix},
\RR^2 b(t,x)=\begin{pmatrix}
18t\rho(t)x_2^2\\
6t
\end{pmatrix},
\end{displaymath}
where $\rho(t)=2+t^3|t|$. Then conditions \eqref{dd1} and  \eqref{ddd1} are satisfied with $N=1$ and $t_1=0$. Moreover, 
\begin{displaymath}
\gamma(t)=R^{-1}(t,x)\RR^2 b(t,x)=\begin{pmatrix}
0\\
\frac{2}{t}
\end{pmatrix}
\end{displaymath}
depends only on $t$. The functions $\gamma_1(t)=0$ and $\gamma_2(t)=\frac{2}{t}$ are analytic in $[-1,1]$ except the point $t_1=0$, where they satisfy condition (i) of Theorem~\ref{th_r}. Since $\gamma_{2,-1}=2$ and all other coefficients $\gamma_{1,i}$ and $\gamma_{2,i}$ vanish, the equation \eqref{eq_th} has the form $k(k-1)-2k=0$. Its roots are $k_1=0$, $k_2=3$. Moreover,  $V_{k,0}=0$ for $k\ge1$, which implies the equality \eqref{eq1_th}. Therefore, all conditions of Theorem~\ref{th_r} hold for the functions $\gamma_1(t)$ and $\gamma_2(t)$. Thus, due to Theorem~\ref{T3}, the system \eqref{sys_ex_last} is locally analytically linearizable in $\R^2$ on the time interval $[-1,1]$. Obviously, in this case $y_1(t)=1$ and $y_2(t)=t^3$ are solutions of the equation \eqref{dif_eq}, which takes the form $\ddot y=\frac{2}{t}\dot y$. Hence, the system \eqref{sys_ex_last} can be transformed to the linear driftless system
\begin{displaymath}
\dot x_1=u, \ \ \dot x_2=t^3u.
\end{displaymath}
\end{ex}

\end{document}